\documentclass[a4paper,12pt,final]{article}
\usepackage{ArXivEN}

\newcommand{\Td}{\mathcal{T}_\Delta} 
\newcommand{\T}{\mathcal{T}_{K}}
\newcommand{\tT}{\widetilde{\mathcal{T}}_{K,\tau}}
\newcommand{\tTin}{\widetilde{\mathcal{T}}_K^{\text{in}}}
\newcommand{\tTn}{\widetilde{\mathcal{T}}_n}
\newcommand{\Tn}{\mathcal{T}_n}
\newcommand{\Tin}{\mathcal{T}_K^{\text{in}}}
\newcommand{\Tle}{\mathcal{T}_{K}^{\text{le}}}
\newcommand{\tTle}{\widetilde{\mathcal{T}}_{K,\tau}^{\text{le}}}
\newcommand{\gd}{g_{\Delta}}
\newcommand{\gtn}{\widetilde{g}_n}
\newcommand{\gtin}{\widetilde{g}_K^\text{in}}
\newcommand{\gtle}{\widetilde{g}_{K,\tau}^\text{le}}
\newcommand{\gt}{\widetilde{g}_{K,\tau}}

\date{}
\title{An Attack on Flexibility and Stoker's Problem}
\author{Maria Hempel}

\hypersetup{
    pdftitle={An Attack on Flexibility and Stoker's Problem},
    pdfauthor={Maria Hempel}
}

\begin{document}
 
\maketitle
\begin{abstract}
In view of solving questions of geometric realizability of polyhedra under given geometric constraints, we parametrize the moduli-space of simply connected polyhedra in terms of their face and their dihedral angles. We recompute its dimension at points of smoothness modulo a combinatorial conjecture.
\end{abstract}

\tableofcontents
\clearpage

\section{Introduction}

\begin{figure}[h]
\centering
\includegraphics{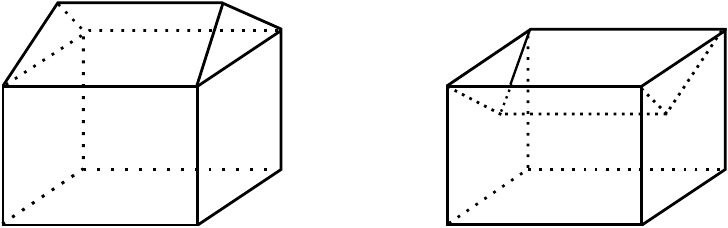}
\caption{Two houses, one with a roof, and one with a roof pool.}
\label{fig:noncongruenthouses}
\end{figure}

\begin{figure}[h] 
\centering
\includegraphics{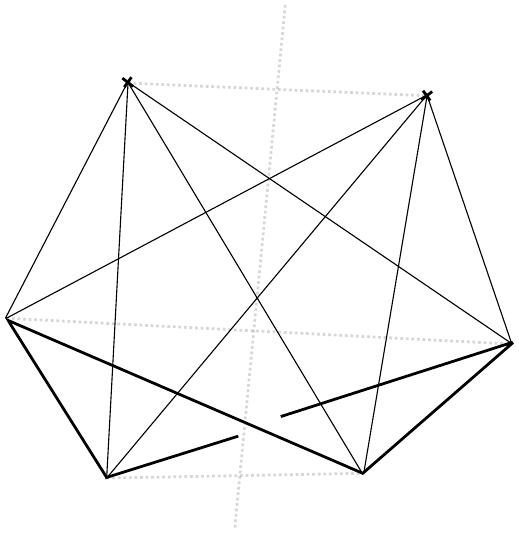}
\caption{Bricard's octahedron of type I intersects itself.}
\label{fig:bricard}
\end{figure}

Polyhedra are simple objects that everybody knows and that are widely used. Nevertheless, fundamental questions about their shape remain open. For instance, we would like to know how the shape of a polyhedron depends on the shape of its faces, and this seemingly simple question still eludes us. It is easy to come up with examples of polyhedra, that show that the shape of the faces does not determine the shape of the polyhedron, see Figure~\ref{fig:noncongruenthouses}.
Looking at the examples in Figure~\ref{fig:noncongruenthouses}, and similar ones, we might be led to think that the only non-uniqueness of the shape of a polyhedron comes about by reflections through planes containing a closed edge-path. In particular, there would be only finitely many shapes of a polyedron with given combinatorics and faces. However, one can even construct continuous families of such, one of which is shown in Figure~\ref{fig:bricard}.

\begin{figure}[t]
\centering
\includegraphics{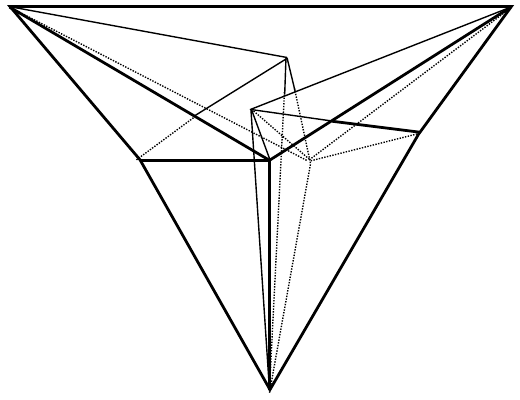}
\caption{Steffen's polyhedron.}
\label{fig:steffen}
\end{figure}

We call such polyhedra \mayhem{flexible}, and the above example along with two other types was the first one found~\cite{bricard:octaedre_articule:1897}. 
With some more maneuvering one can also produce examples of flexible polyhedra without self-intersections, which was done by Connelly~\cite{connelly:an_immersed_flexie:1976,connelly:a_counterexample:1977} in the late seventies. In Figure~\ref{fig:steffen}, Steffen's polyhedron consists of two Bricard octahdera with two removed faces and joined by a flap. A printable folding pattern can be found on cutupfoldup.com. Beautiful flexible polyhedra that are variations on the theme of those of Bricard were done by Nelson~\cite{nelson:generalizations_of_bricard_octahedra,nelson:extending_bricard_octahedra}. 

The question whether a given polyhedron is flexible is asking, whether there is a continuous family of realizations for given fixed face angles. Similarly, one might ask whether the dihedral angles of a polyhedron determine its shape, and again we know very little. We do know that infinitesimally the shape cannot change when the dihedrals are fixed~\cite{mazzeo_and_montcouquiol:Infinitesimal_rigidity_of_cone-manifolds_and_the_Stoker_Problem}. We also know that the shape is unique for certain restrictive subclasses of convex polyhedra with fixed dihedrals~\cite{karcher:Remarks_on_polyhedra_with_given_dihedral_angles,milka:On_a_certain_conjecture_of_Stoker,pogorelov:On_a_problem_of_Stoker}. The question whether this is true for all convex polyhedra is known as Stoker's Conjecture and is still open.

Indeed, one is often led to ask for any set of given geometric data, such as subsets of face angles, dihedrals or edgelengths, what the space of possible realizations is. In the following we attack these questions by giving a new description of the space of shapes that is natural to the problem. We note that spaces of shapes of polyhedra have been studied, see~\cite{richter-gebert:realization_spaces_of_polytopes:1996}, and references therein, and Whiteley~\cite{whiteley:how_to_describe_or_design_a_polyhedron}. Here, the space of shapes is described by studying the possible liftings to three-space of planar graphs that may be considered projections of polyhedra, and several conlcusions about which geometric data can be determined from which other data have been drawn from it. In~\cite{thurston:shapes_of_polyhedra:1998} Thurston Sr. studies the space of flat-cone surfaces with given convex cone-singularities, each of which by a theorem of Alexandrov~\cite{alexandrov:convex_polyhedra} and by Cauchy's Rigidity Theorem~\cite{cauchy:sur_les_polygones_et_polyedres:1813,Hadamard:Erreurs_des_mathematiciens,lebesgue:demonstration_complete_du_theoreme_de_Cauchy} can be uniquely realized as a convex polyhedron. However, this proof is non-constructive, so that it does not help to relate geometric contraints. In the following we describe the shape of a polyhedron in terms of its face angles and its dihedral angles by describing it as a flat-cone surface with specified triangulation and by describing the polyhedral cones defined around each vertex.

By a polyhedron we really mean the following.

\begin{definition}
	A \mayhem{combinatoric} is a nonempty, consistently oriented abstract simplicial two-complex, whose geometric realization is homeomorphic to a closed connected surface.
\end{definition}

\begin{definition}
	A \mayhem{polyhedron} is a simplexwise linear and simplexwise injective map $P\colon\abs{K}\to \R^3$ from the geometric realization $\abs{K}$ of a combinatoric $K$ into euclidean $\R^3$.
\end{definition}

We call the images of a polyhedron of its zero-, one- and two-simplices its \mayhem{vertices}, \mayhem{edges}, and \mayhem{faces}, respectively. Their numbers will always be denoted by the letters $V$, $E$ and $F$, respectively. Note that our polyhedra are always triangulated and oriented. Furthermore, a polyhedron may intersect itself and even degenerate into a plane, but its faces and edges are non-degenerate. We call two faces or edges \mayhem{adjacent}, whenever their preimages are adjacent as simplices. A polyhedron $P\colon \abs{K}\rightarrow \R^3$ defines two special sets of angles, namely

\[
	\delta_P \in [0,2\pi)^{K_1} 
	\quad\quad\text{ and }\quad\quad 
	\sigma_P\in (0,\pi)^{C_K},
\]where $\delta_K$ are the dihedral angles between adjacent faces of $P$, labeled by the one-simplices of $K$, and $\sigma_K$ are the angles of the faces of $P$ labeled by the set of \mayhem{corners of $K$}
\[
	C_K :=\Set{ (ijk) \mid (ijk)\text{ is an order of }\set{ijk}\in K_2}
	\Big/\sim,
	\]
	where $\sim$ identifies an order $(ijk)$ with the order $(kji)$.

\begin{figure}[t]
\centering
\includegraphics{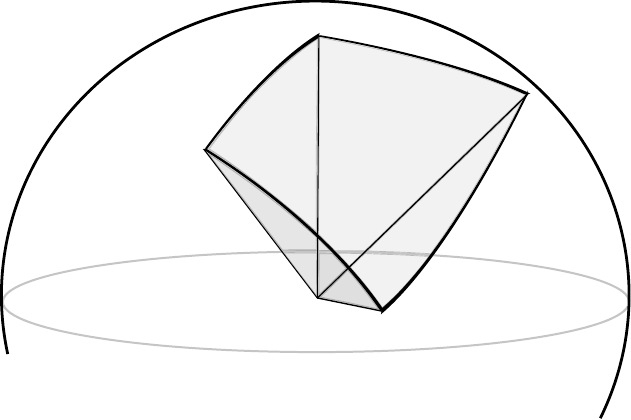}
\caption{Legendre's trick}
\label{fig:cauchystrick}
\end{figure}

It is easy to prove that the angles $\delta_P$ and $\sigma_P$ completely determine the shape of a polyhedron, so that they are natural coordinates for a space of shapes. To define such a space we would like to know which vectors $\sigma\in (0,\pi)^{C_K}$ and $\delta \in [0,2\pi)^{K_1}$ are the surface and dihedral angles of a polyhedron. To answer this question we observe that every polyhedron defines a triangulated flat-cone surface. In Section~\ref{sec:Parametrizing_intrinsic_constraints}, we parametrize the space $\Tin$ of flat-cone surfaces that are triangulated by a combinatoric $K$ in terms of its surface angles, by parametrizing the space of euclidean triangles for each of the faces. In addition, we note that every polyhedron defines a polyhedral cone around each of its vertices. In Section~\ref{subsec:Geometric_realizations_of_a_polyhedral_cone} we parametrize the space $\Tn$ of polyhedral cones with $n$ edges in terms of its surface and its dihedral angles at smooth points. We note that Kapovich and Millson~\cite{kapovich_millson:hodge_theory_and_the_art_of_paper_folding} have allready studied the space of polyhdral cones using Hodge Theory. We use a trick invented by Legendre, and intersect the polyhedral cone with a unit sphere centered at its vertex, see Figure~\ref{fig:cauchystrick}. This intersection is a spherical polygon whose edgelengths are the surface angles of the polyhedral cone, and whose angles are the dihedral angles of the polyhedral cone. We study the space of spherical $n$-gons by triangulating them and parametrizing the space of spherical triangles~\ref{subsec:Geometric_realizations_of_a_spherical_triangle}. We amass the spaces of polyhedral cones defined for each zero-simplex of a combinatoric $K$ and identify dihedrals according to $K$ into a space $\Tle$ in Section~\ref{sec:The_collection_of_polyhedral_cones}. In summary, we observe that the surface and dihedral angles $\sigma_P$ and $\delta_P$ of a polyhedron $P$ are the surface angles of a flat-cone surface with triangulation $K$ and the surface and dihedral angles of a certain collection of polyhedral cones that share dihedral angles at the one-simplices of $K$, i.e., that 
\[
	(\sigma_P,\delta_P) \in \T := (\Tin \times [0, 2\pi)^{K_1}) \cap \Tle.
\]

We find and prove in Section~\ref{sec:Wrapping_it_up} that this is also sufficient for a set of angles to arise as those of a simply connected polyhedron $P$.
\begin{theorem}[Main]
If $K$ is a combinatoric such that $\abs{K}$ is simply connected, then 
\[
	(\sigma,\delta) \in (0,\pi)^{C_K} \times [0,2\pi)^{K_1}
\]
are the surface and dihedral angles $\sigma = \sigma_P$ and $\delta = \delta_P$ of a polyhedron $P$ if and only if $(\sigma,\delta) \in \T$.
\end{theorem}

Furthermore, we can compute some of the following dimensions.

\begin{theorem} 
At points of smoothness the spaces $\Tin$, $\Tn$, $\Tle$ and $\T$ are real analytic manifolds of the following dimensions
\[
	\dim \Tn = 2n-3, \quad \text{ and } \quad \dim \Tle = 2E + 6g - 6,
\]
where $g$ is the genus of $\abs{K}$, and it is well known and we recompute modulo a well decorated conjecture~\ref{thm:elimination}, that
\[
	\dim \T = \dim \Tin  = E-1.
\]
\end{theorem}

\section{Parametrizing intrinsic constraints}\label{sec:Parametrizing_intrinsic_constraints}

Intrinsically, a polyhedron is a \mayhem{flat-cone surface}, that is, with the induced metric from the ambient euclidean three-space, it is a surface that is locally isometric to the euclidean plane, with the exception of neighborhoods of finitely many points that are isometric to a euclidean cone. Furthermore, the edges of $P\colon \abs{K}\rightarrow \R^3$ define a triangulation of its flat-cone surface $S_P$, whose vertices are the cone points of $S_P$. We note that $P$ defines a map $\ell_P \in \R_{+}^{K_1}$ of edgelengths of $P$, from which $S_P$ can be entirely reconstructed. In other words, the faces of a polyhedron define a collection of euclidean triangles that share edges as prescribed by $K$, and it is sufficient to parametrize the space of euclidean triangles for each of the triangles of $\abs{K}$, to describe the possible intrinsic geometries of $P$.

\subsection{Geometric realizations of a euclidean triangle}\label{subsec:geometric_realizations_of_a_euclidean_triangle}

We parametrize the geometric realizations of a nondegenerate euclidean triangle in terms of its edgelengths and its angles as follows.

\begin{definition}
Let
\[
	\gd\colon \R_+^3\times (0,\pi)^3 \longrightarrow \R^3
\]
be defined by
\[
	g_\Delta(a,b,c,\alpha,\beta,\delta) 
	= \begin{pmatrix}
	a\cos\beta + b\cos\alpha - c\\
	b\sin\alpha - a\sin\beta\\
	\alpha + \beta + \gamma -\pi\end{pmatrix},
\]
and set 
\[
	\Td:= g_\Delta^{-1} (0).
\]
\end{definition}

\begin{proposition}\label{thm:triangle} 
$\mathcal{T}_\Delta$ is an analytic submanifold of $\R_+^3\times (0,\pi)^3$ of dimension~$3$, and $(a,b,c,\alpha,\beta,\gamma)\in\mathcal{T}_\Delta$ if and only if $a$, $b$ and $c$ are the edgelengths of a nondegenerate euclidean triangle with opposite angles $\alpha$, $\beta$ and $\gamma$, respectively.
\end{proposition}

\begin{proof}

We compute
\[
\D g_\Delta = \begin{pmatrix}
		\cos\beta & \cos\alpha & -1 & -b\sin\alpha & -a\sin\beta &0\\
		-\sin\beta& \sin\alpha &  0 &  b\cos\alpha & -a\cos\beta &0\\
		0	  & 0	       &  0 &  1           &  1		 &1
		\end{pmatrix}
\]
and find that, for example, columns $3$, $2$ and $6$ form an upper-triangular $3\times 3$ submatrix whose diagonal entries are nonzero on $\R_+^3\times (0,\pi)^3$ and in particular on $\mathcal{T}_\Delta$. Hence, on $\Td$,
\[
\rk \D g_\Delta =3
\]
on $\mathcal{T}_\Delta$. Therefore, by the Preimage Theorem, confer for instance~\cite{guillemin_pollack:differential_topology:2010}, $\mathcal{T}_\Delta$ is an analytic submanifold of $\R_+^3\times (0,\pi)^3$ of dimension
\[
\dim\mathcal{T}_\Delta = \dim(\R_+^3\times (0,\pi)^3) -\rk \D g_\Delta=3.
\]

	Note that if $a$, $b$, $c$ and $\alpha$, $\beta$ and $\gamma$ are the edgelengths and opposite angles of a euclidean triangle then it is well known that they satisfy the claimed equations. The converse is an easy exercise in plane euclidean geometry.
\end{proof}

\subsection{Geometric realizations of triangulated flat cone surfaces}\label{subsec:geometric_realizations_of_flat_cone_surfaces_with_prescribed_triangulation}

Let $K$ be a combinatoric and for every map $\ell\in\R_+^{K_1}$ and map $\sigma\in (0,\pi)^{C_K}$ abbreviate 
\[
	\ell_{ij}:=\ell(\set{ij}) 
	\quad \text{ and } \quad
	 \sigma_{ijk} := \sigma\left(\left( ijk \right)  \right).
\] For every two-simplex $\set{ijk}\in K_2$ we parametrize a euclidean triangle with edgelengths $\ell_{ij}$, $\ell_{jk}$ and $\ell_{ik}$ and opposite angles $\sigma_{ikj}$, $\sigma_{jik}$ and $\sigma_{ijk}$, respectively, as follows.

\begin{definition}
Let
\[
	\gtin: \R_+^{K_1} \times (0,\pi)^{C_K}  \longrightarrow (\R^3)^{K_2}
	\] be the function defined by
\[
	\left(\gtin(\ell,\sigma)\right)(\set{ijk}) := \gd(\ell_{ij},\ell_{jk},\ell_{ki},\sigma_{ikj},\sigma_{jik},\sigma_{kji}) 
	\]
and let
\[
	\tTin := \left( \widetilde{g}_K^\text{in} \right)^{-1}(0).
	\]
\end{definition}

\begin{remark}
	Note that $\left(\gtin\right)(\set{ijk})$ is neither symmetric in $\ell_{ij}$, $\ell_{jk}$ and $\ell_{ki}$, nor in $\sigma_{ijk}$, $\sigma_{jki}$ and $\sigma_{kij}$. This means that strictly speaking $\widetilde{g}_K^\text{in}$ depends on the choice of a cyclic order of the zero-simplices of every two-simplex. However, any solution$(a,b,c,\alpha,\beta,\gamma)\in\mathcal{T}_\Delta$, also defines a solution $(b,c,a,\beta,\gamma,\alpha)\in\mathcal{T}_\Delta$.  Hence, even though the function $\gtin$ depends on the choice of an order of the zero-simplices of every two-simplex of $K$, the set of solutions $\tTin$ does not, and we omit the choice in the notation.
\end{remark}

\begin{proposition}\label{thm:intrinsic}
	Let $K$ be a combinatoric. Then $\tTin$ is an analytic submanifold of $(0,\pi)^{C_K}\times\R_+^{K_1}$ of dimension $E$, and $(\ell,\sigma)\in\tTin$ if and only if $\ell$ are the edgelengths of a flat cone surface triangulated by $K$ and $\sigma$ are its angles.
\end{proposition}

\begin{proof}
	Each two-simplex $\set{ijk}\in K_2$ contributes three rows to $\D \gtin$, which we denote by $r_{ijk}$, $s_{ijk}$ and $t_{ijk}$ and which, omitting the columns that are zero, look as follows

\[ 	\begin{pmatrix}
		r_{ijk}\\
		s_{ijk}\\
		t_{ijk}
	\end{pmatrix}
	:=\begin{pmatrix}
	-\ell_{ij}\sin\sigma_{ijk} & -\ell_{ik}\sin\sigma_{jki} & 0 & \cos\sigma_{ijk} & \cos\sigma_{jki} & -1\\
	 \ell_{ij}\cos\sigma_{ijk} &  -\ell_{ik}\cos\sigma_{jki} & 0 & \sin\sigma_{ijk} & -\sin\sigma_{jki}&  0\\
		 1	     &  1	   & 1 & 0         & 0          &  0
	\end{pmatrix},
	\]
	where the first three columns shown correspond to the partial derivatives in the angles of $\set{ijk}$. Since $\gtin$ only depends on the angles of $\set{ijk}$ on $\set{ijk}$, the other entries of the first three columns are zero. Suppose now for every $\set{ijk}\in K_2$ there were constants $\lambda_{ijk}, \mu_{ijk}, \nu_{ijk}\in\R$, such that
\begin{align}
	\sum\lambda_{ijk} r_{ijk} + \sum\mu_{ijk} s_{ijk} + \sum\nu_{ijk} t_{ijk} &= 0.\label{eq:lambdamunu} 
\end{align}
Then, since each $t_{ijk}$ contains an entry $1$ in a column whose other entries are all $0$, it follows that $\nu_{ijk} = 0$ for all $\set{ijk}\in K_2$. Together with~\ref{eq:lambdamunu} this implies 
\begin{align}
	\sum\lambda_{ijk} r_{ijk} + \sum\mu_{ijk} s_{ijk} &= 0.\label{eq:lambdamu} 
\end{align}
Therefore, reading~\ref{eq:lambdamu} for the first column means that
\begin{align}
 -\lambda_{ijk} \ell_{ij}\sin\sigma_{ijk} + \mu_{ijk} \ell_{ij}\cos\sigma_{ijk} = 0,\label{eq:1stcol}
\end{align}
 and reading it for the second column means that
\begin{align}
 -\lambda_{ijk} \ell_{ik}\sin\sigma_{jki} - \mu_{ijk} \ell_{ik}\cos\sigma_{jki} = 0.\label{eq:2ndcol}
\end{align}
Since $\ell_{ij}$ and $\ell_{ik}$ are nonzero, we can divide Equation~\ref{eq:1stcol} and Equation~\ref{eq:2ndcol} by $\ell_{ij}$, respectively $\ell_{ik}$, and get
\begin{align}
	\mu_{ijk}\cos\sigma_{ijk} &= \lambda_{ijk}\sin\sigma_{ijk},\label{eq:mucos}
\end{align}
respectively
\begin{align}
	\lambda_{ijk}\sin\sigma_{jki} &= -\mu_{ijk}\cos\sigma_{jki}.\label{eq:lamsin}
\end{align}
If $\cos\sigma_{ijk} = 0$, then, since $\sigma_{ijk}\in(0,\pi)$ it follows that $\sigma_{ijk} = \tfrac{\pi}{2}$ and~\ref{eq:mucos} becomes
\[
	0 = \lambda_{ijk}.
	\]
If $\cos\sigma_{ijk} \neq 0$, then we can divide Equation~\ref{eq:mucos} by $\cos\sigma_{ijk}$ and get
\begin{align}
	\mu_{ijk} &= \lambda_{ijk}\tan\sigma_{ijk}.\label{eq:mu}
\end{align}
Similarly, we note that $\sin\sigma_{jki} \neq 0$ for $\sigma_{jki}\in (0,\pi)$, so that Equation~\ref{eq:lamsin} becomes
\begin{align}
	\lambda_{ijk} &= - \mu_{ijk} \frac{1}{\tan\sigma_{jki}}.\label{eq:lam}
\end{align}
Therefore, substituting$\lambda_{ijk}$ of Equation~\ref{eq:lam} into Equation~\ref{eq:mu}, we get
\begin{align}
	\mu_{ijk} &= -\mu_{ijk}\frac{\tan\sigma_{ijk}}{\tan\sigma_{jki}}.\label{eq:mutan}
\end{align}
If $\mu_{ijk}$ were nonzero, then we could divide Equation~\ref{eq:mutan} by it and would get
\begin{align}
	\tan\sigma_{jki} &= -\tan\sigma_{ijk},\label{eq:tansig}
\end{align}
which for $\sigma_{ijk},\sigma_{jki}\in (0,\pi)$ means
\begin{align}
	\pi -\sigma_{kji} = \sigma_{ijk}.\label{eq:twoanglessumpi}
\end{align}
But, on $\tTin$ we have that $\sigma_{ijk}$, $\sigma_{jki}$ and $\sigma_{kij}$ add up to $\pi$. Therefore,~\ref{eq:twoanglessumpi} implies that $\sigma_{kij}=0$, which contradicts $\sigma_{kij}\in (0,\pi)$. It follows that $\mu_{ijk} = 0$.

	Similarly, substituting $\mu_{ijk}$ of Equation~\ref{eq:mu} into Equation~\ref{eq:lam}, we get
\begin{align}
	\lambda_{ijk} &= -\lambda_{ijk}\frac{\tan\sigma_{ijk}}{\tan\sigma_{jki}}, \label{eq:lamtan}
\end{align}
so that if $\lambda_{ijk}$ were nonzero, then~\ref{eq:lamtan} would also imply~\ref{eq:tansig}, which again is impossible. It follows that $\lambda_{ijk}=0$.

Hence, we have shown that on $\widetilde{\mathcal{T}}_K^\text{in}$ the rows of $\D \gtin$ are linearly independent and that it is therefore of full rank. This proves that $0$ is a regular value of $\gtin$ so that by the preimage theorem, confer for instance~\cite{guillemin_pollack:differential_topology:2010}, we get that $\tTin$ is a smooth submanifold of dimension
\begin{align*}
	\dim\tTin
	&= \dim \left((0,\pi)^{C_K}\times \R_+^{K_1}\right)-\dim \left((\R^3)^{K_2}\right)\\
	&= (3F + E) - 3F = E.
\end{align*}
The rest of the claim follows by construction of $\gtin$ and Proposition~\ref{thm:triangle}.
\end{proof}

It is intuitively clear that each of the edgelengths of a flat-cone surface can be chosen independently, which we have shown algebraically in the previous proposition. Similarly, we show below that the set of angles of a triangulated flat-cone surface together with one edgelength determines all of its edgelengths. We note that this means in particular that there are more constraints on the angles of a flat-cone surface, than that they sum to $2\pi$ on each face individually. Indeed, since they share edges the angles of two faces are related by a sequence of laws of the sines for a path of faces between them.  Below we use these constraints to construct a diffeomorphism that incorporates all the information needed to fully describe the set of angles of a triangulated flat-cone surface.

\begin{definition}
Let $p_\Sigma$ be the projection of $\R_{+}^{K_1}\times (0,\pi)^{C_K}$ onto $(0,\pi)^{C_K}$ and set
\[
	 \Tin := p_\Sigma \left(\tTin\right).
	\]
\end{definition}

\begin{proposition}\label{thm:intrinsicdiffeo}
For every $\set{ij}\in K_1$ there is a diffeomorphism 
\[
	\varphi_{ij}\colon \mathcal{T}_K^\text{in} \times \R_+ \longrightarrow \widetilde{\mathcal{T}}_K^\text{in}.
	\]
	In particular, $\mathcal{T}_K^\text{in}$ is a smooth submanifold of $(0,\pi)^{C_K}$ of dimension
\[
	\dim \mathcal{T}_K^\text{in} = E-1.
	\]
\end{proposition}

\begin{proof}
	We start with the following observation for a given $(\sigma,x)\in\tTin$. For every $\set{lm}\in K_1\setminus \set{ij}$ choose a path $\gamma_{lm}$ of one-simplices from $\set{ij}$ to $\set{lm}$, that is, a finite sequence of consecutively adjacent one-simplices, whose first element is $\set{ij}$ and whose last element is $\set{lm}$. Denote this choice of paths by $\gamma := \left(\gamma_{lm}\right)_{\set{lm}\in K_1\setminus\set{ij}}$.
	\begin{figure}[tpb]
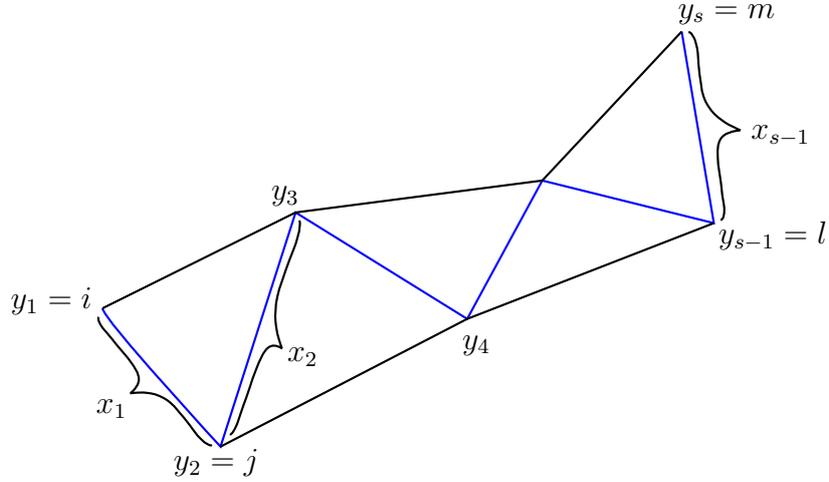

		\begin{center}
			\op{sinealong2}{
				\put(0,20){$y_1 = i$}
				\put(20,0){$y_2 = j$}
				\put(32,33.5){$y_3$}
				\put(55.5,15){$y_4$}
				\put(87,28){$y_{s-1}=l$}
				\put(82,56){$y_s = m$}
				\put(10.5,7){$x_1$}
				\put(34,13.5){$x_2$}
				\put(91,41){$x_{s-1}$}
					}
		\end{center}
		\caption{In blue: a path $\gamma_{lm}$.}
		\label{fig:sinealong}
	\end{figure}
Denote now
\[
	y_1=i, y_2=j, y_3,\ldots, y_{s-1}=l, y_s= m 
	\]
the zero-simplices transversed by $\gamma_{lm}$ from $i$ to $m$. Note that $y_1$,\ldots, $y_s$ depend on $\gamma$ and on $\set{lm}$ but we omit this dependence in the notation for the sake of readability. Denote 
\[
	x_k:=x_{y_k y_{k+1}},
	\]
and compare with Figure~\ref{fig:sinealong}.
	By Proposition~\ref{thm:intrinsic} for every two-simplex $\set{i'j'k'}\in K_2$ the values $x_{i'j'}$, $x_{j'k'}$ and $x_{k'i'}$ are the edgelengths of a euclidean triangle with opposite angles $\sigma(j'k'i')$, $\sigma(k'i'j')$ and $\sigma(i'j'k')$. Therefore, they obey the law of the sines. In particular, we have,
	\[
		x_{k+1} = \frac{\sin\sigma(y_{k+1}y_{k-1}y_{k+2})}{\sin\sigma(y_{k}y_{k+2}y_{k+1})} x_k,
		\]
for $k\in\set{1, 2, \ldots,s-2}$, so that
\begin{align}
	x_{lm}=x_{s-1}  =\frac{\sin\sigma(y_{s-1}y_{s-2}y_s)}{\sin\sigma(y_{s-2}y_sy_{s-1})} \cdots \frac{\sin\sigma(y_2y_1y_3)}{\sin\sigma(y_1y_3y_2)} x_{ij}.\label{eq:sinealong}
\end{align}

Now, note that $\sigma\in\Tin$ means that $\sigma$ are the angles of a flat cone surface triangulated by $K$, and in particular, that there exists a flat cone surface $s\in\R_+^{K_1}$ triangulated by $K$, whose angles are $\sigma$. Note, furthermore, that if $s$ is a flat cone surface triangulated by $K$ with angles $\sigma$, then for every $\lambda\in\R_+$ the tuple $(\sigma,\lambda s)$ also satisfies the defining equations of $\tTin$, so that we have that $(\sigma,\lambda s)\in\tTin$. Therefore, for every $(\sigma,\ell)\in\Tin\times\R_+$ there is a flat cone surface triangulated by $K$, whose edgelength of $\set{ij}$ is $\ell$. By our previous observation its other edgelengths are given by 
\begin{align}
	d_{lm}^\gamma(\sigma, \ell):= \frac{\sin\sigma(y_{s-1}y_{s-2}y_s)}{\sin\sigma(y_{s-2}y_sy_{s-1})} \cdots \frac{\sin\sigma(y_2y_1y_3)}{\sin\sigma(y_1y_3y_2)} \ell,\label{eq:sinealongell}
\end{align}
for $\set{lm}\in K_1\setminus\set{ij}$ and by $\ell$ for $\set{ij}$.
We have now seen, that the map
\[
	\varphi_{ij}^\gamma\colon \mathcal{T}_K^\text{in}\times \R_+
	\longrightarrow 
	\widetilde{\mathcal{T}}_K^\text{in}
	\]
given by
\[
	\varphi_{ij}^\gamma \left(\sigma, \ell\right) := (\sigma, d^\gamma(\sigma, \ell)),
	\]
	where $d^\gamma(\sigma,\ell)$ sends each one-simplex $\set{lm}\in K_1\setminus\set{ij}$ to $d_{lm}^\gamma(\sigma,\ell)$ and $\set{ij}$ to $\ell$, is well defined.

For $\sigma\in (0,\pi)^{C_K}$ the denominators in~\ref{eq:sinealongell} do not vanish and therefore $\varphi_{ij}^\gamma$ is smooth. Furthermore, we find that
\[
	(\varphi_{ij}^\gamma)^{-1}(\sigma,x) := (p_{C_K}(\sigma,x), x_{ij})
	\]
	defines a smooth inverse function for $\varphi_{ij}^\gamma$ and that it is thus a diffeomorphism. The inverse function does not depend on the choice of  $\gamma$. It follows that $\varphi_{ij}^\gamma$ does not depend on the choice of $\gamma$ and we can henceforth drop the reference to it.
\end{proof}

\section{Parametrizing local extrinsic constraints}\label{sec:Parametrizing_local_extrinsic_constraints}

The intersection of a possibly rescaled polyhedron $P$ with a unit sphere centered at a vertex is a spherical polygon whose edge-lengths are the angles of $P$ around that vertex, and whose angles are the dihedral angles of $P$ around that vertex. We use this observation to relate the surface and dihedral angles explicity in the following, by triangulating each of the spherical polygons thus defined, parametrizing the spherical triangles of the triangulation, and then eliminating the auxiliary variables.

\subsection{Geometric realizations of a spherical triangle}\label{subsec:Geometric_realizations_of_a_spherical_triangle}

We consider the unit sphere $S^2$ in $\R^3$ centered at the origin. For an introduction to spherical geometry see~\cite{berger:geometry2}. Recall that the intersection of $S^2$ with a plane passing through the origin is called a (spherical) \mayhem{line}. A line on $S^2$ is a curve, so that it also makes sense to speak of line segments.

\begin{definition}
A \mayhem{circuit} is a simplicial complex $Z=Z_0\cup Z_1$ with zero-simplices

\[	Z_0 = \Set{\set{x_1}, \dots, \set{x_n}},
	\]
where $n>2$ and with one-simplices
\[
	Z_1 = \Set{\set{x_i x_{i+1}} \mid i\in\set{1,\ldots,n-1}}\cup\Set{\set{x_n x_1}}. 
	\]
When we want to bring out the number $n$ of zero-simplices of $Z$ we also call $Z$ an \mayhem{$n$-circuit}.
\end{definition}

\begin{definition}
	A \mayhem{spherical polygon} is a map $p\colon \abs{Z}\to S^2$ from the geometric realization of a circuit $Z$ into the unit sphere $S^2$, that maps one-simplices to line segments of lengths contained in $(0,\pi)$. If $Z$ is an $n$-circuit, then we also call $p$ a \mayhem{spherical $n$-gon}.
\end{definition}

We call the images of the zero-simplices of a spherical $n$-gon its \mayhem{vertices} and the images of its one-simplices its \mayhem{edges}. It is easy to see that a spherical $3$-gon whose vertices define linearly independent vectors in $\R^3$ is injective. In particular, such a $3$-gon is a simple closed curve on $S^2$ and by the Jordan Curve Theorem, divides the unit sphere into two regions. 

\begin{definition}
The union of a spherical $3$-gon whose vertices define linearly independent vectors in $\R^3$ with one of the two regions that it divides $S^2$ into is called a \mayhem{spherical triangle}.
\end{definition}

By definition a spherical triangle has edgelengths contained in $(0,\pi)$. Its \mayhem{angle} at a vertex is the angle of the tangent vectors to the two adjacent edges that projects onto the interior of the spherical triangle. We easily find that the angles of a spherical triangle defined as above are either all contained in $(0,\pi)$, or all contained in $(\pi,2\pi)$. This together with the spherical cosine rule motivates the following.

\begin{definition}
	Let
\[
	g_3 \colon (0,\pi)^6\cup \left((0,\pi)^3\times(\pi,2\pi)^3\right) \longrightarrow \R^3
	\]
be defined by
\[ 
	g_3(a,b,c,\alpha,\beta,\gamma) 
	:=\begin{pmatrix}
		\cos b\cos c +\sin b \sin c\cos\alpha -\cos a	\\
		\cos c\cos a +\sin c \sin a\cos\beta  -\cos b	\\
		\cos a\cos b +\sin a \sin b\cos\gamma -\cos c	
	\end{pmatrix},
	\]
and set
\[
	\mathcal{T}_3 := g_3^{-1}(0).
\]
\end{definition}

\begin{proposition}\label{thm:sphericaltriangle}
The set $\mathcal{T}_3$ is an analytic manifold of dimension $3$ and $(a,b,c, \alpha,\beta,\gamma) \in \mathcal{T}_3$ if and only if there is a spherical triangle with edgelengths $a$, $b$ and $c$ and opposite angles $\alpha$, $\beta$ and $\gamma$, respectively.
\end{proposition}

\begin{proof}
One easily finds a submatrix of $\D g_3$ which is of full rank on the domain of $g_3$. Indeed, any triple of partial derivatives by variables occuring together in a spherical congruence theorem, is of full rank $3$. The first part of the claim then follows from the Preimage Theorem. The second part of the claim is an easy exercise in spherical geometry.
\end{proof}

\subsection{Geometric realizations of a polyhedral cone in $\R^3$}\label{subsec:Geometric_realizations_of_a_polyhedral_cone}

Let $K$ be a combinatoric, and $k\in K_0$ be a zero-simplex. Then we call the smallest simplical two-subcomplex of $K$ that contains all the simplices of $K$ that conain $k$, the \mayhem{abstract cone} of $k$, and denote it by $C(k)$. The image of a simplex-wise linear and injective map $\abs{C(k)}\longrightarrow \R^3$ defines a union of euclidean triangles in $\R^3$ with a common vertex, and the union of the planar sectors defined by them is called a \mayhem{polyhedral cone}.
 In this way a polyhedron $P\colon \abs{K}\longrightarrow \R^3$ defines a polyhedral cone $C(v)$ for each of its vertices $v$. If $k\in K_0$ has valency $n_k$, then we say that the polyhedral cone defined by $P$ for $k$ is an \mayhem{$n_k$-cone}. We call $v$ the vertex of $C(v)$, the half-lines at $v$ defined by the edges of $P$, the \mayhem{edges} of $C(v)$, and the plane sectors defined by the faces of $P$ containing $v$, the \mayhem{faces} of $C(v)$ .

\begin{figure}[tpb]
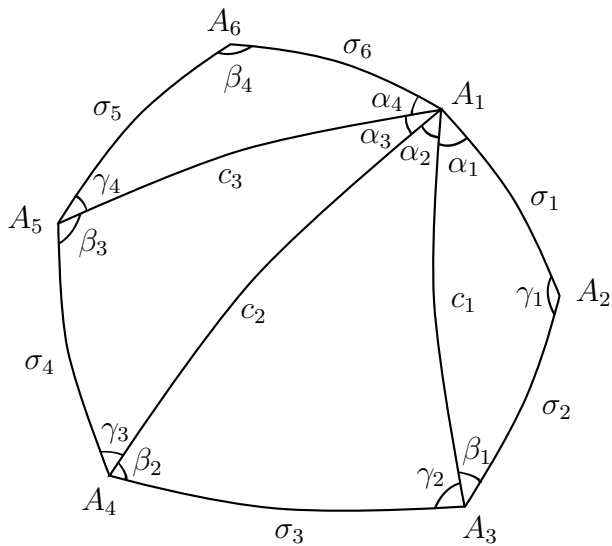

	\begin{center}
		\op{pv}{
			\put(75,78){$A_1$}
			\put(97,43){$A_2$}
			\put(77,1){$A_3$}
			\put(10,6){$A_4$}
			\put(-3,56){$A_5$}
			\put(32,91){$A_6$}
			\put(89,60){$\sigma_1$}
			\put(91,23){$\sigma_2$}
			\put(44,1){$\sigma_3$}
			\put(0,31){$\sigma_4$}
			\put(12,76){$\sigma_5$}
			\put(56,87){$\sigma_6$}
			\put(75,42){$c_1$}
			\put(38,40){$c_2$}
			\put(34,64){$c_3$}
			\put(74.5,66){$\alpha_1$}
			\put(66,68){$\alpha_2$}
			\put(59,71){$\alpha_3$}
			\put(61,77.5){$\alpha_4$}
			\put(86.5,43){$\gamma_1$}
			\put(69,11){$\gamma_2$}
			\put(13.5,19){$\gamma_3$}
			\put(12,63){$\gamma_4$}
			\put(77,15){$\beta_1$}
			\put(19,13){$\beta_2$}
			\put(9.5,52){$\beta_3$}
			\put(35,81){$\beta_4$}
			}
	\end{center}
	\caption{A triangulated spherical polygon}
	\label{fig:pk}
\end{figure}

The intersection of a polyhedral cone $C(v)$ defined by a polyhedron $P$ at its vertex $v$ of valency $n$ with a small unit sphere is a spherical n-gon, as defined earlier. To give it a triangulation we pick a cyclic labeling $A_1, A_2,\ldots, A_n$ of its vertices and consider the spherical line-segments between $A_1$ and $A_i$ for $i=3,\ldots, n-1$, compare Figure~\ref{fig:pk}. Whenever no two edges of $C(v)$ are colinear then there is a unique shortest spherical line segment from $A_1$ to $A_i$ and its length $c_i$ is contained in $(0,\pi)$. Hence, in this situation the vertices $A_1$, $A_i$ and $A_{i+1}$ for $i=2,\ldots,n-2$ define a spherical three-gon $\Delta_{i-1}$. If in addition no three edges of $C(v)$ are coplanar, then these spherical three-gons are injective and by Jordan's Curve Theorem define two disjoint regions of $S^2$ and the choice of any one of them makes $\Delta_i$ into a spherical triangle.

\begin{definition}
A polyhedral cone is called \mayhem{in general position} if any triple of its edges is linearly independent.
\end{definition}

Therefore, whenever $C(v)$ is in general position, we pick as interior of $\Delta_i$ the region of $S^2$ indicated by the inward pointing normal of $P$, and have defined a \mayhem{triangulation} of $C(v)$. Let $\alpha_i$, $\beta_i$, and $\gamma_i$ be the angles of $\Delta_i$, as shown in Figure~\ref{fig:pk}. 

\begin{definition}
	Let
\[
	\gtn\colon\Omega_3^{n-2}\times \R_+^{K_1}\longrightarrow  \R^{3(n-2)+ n}
\]
be defined by
\[
	\gtn(\alpha,\beta,\gamma,c,\sigma,\delta)
	:=\begin{pmatrix}
		g_3(\sigma_2,\sigma_1, c_1,\alpha_1,\beta_1,\gamma_1)\\
		g_3(\sigma_3,c_1, c_2,\alpha_2,\beta_2,\gamma_2)\\
		\vdots\\
		g_3(\sigma_{n-2},c_{n-4}, c_{n-3},\alpha_{n-3},\beta_{n-3},\gamma_{n-3})\\
		g_3(\sigma_{n-1},c_{n-3},\sigma_n,\alpha_{n-2}\beta_{n-2}\gamma_{n-2})\\
		\delta_1-\left(\alpha_1 + \cdots + \alpha_{n-2}\right)\\
		\delta_2-\gamma_1\\
		\delta_{3}-\beta_{1} -\gamma_{2} \\
		\vdots\\
	        \delta_{n-1}-\beta_{n-3} -\gamma_{n-2} \\
		\delta_n - \beta_{n-2}
	\end{pmatrix}
	\]
and set
\[
	\tTn : = \gtn^{-1}(0)
	\quad\text{ and }\quad
	\Tn : = p_n(\tTn),
\]
where $p_n$ is the projection onto the $\sigma$- and $\delta$-factors of $\Omega_n:=\Omega_3^{n-2}\times (0,2\pi)^{K_1}$.
\end{definition}

\begin{remark} To write $\Omega_3^{n-2}\times \R_+^{K_1}$, as we do above, as the domain for $\gtn$ is an abuse of notation, since $\gtn$ draws its variables from its domain as written in a different order. For a long clean ennummeration we refer the reader to the thesis version of this paper.
\end{remark}

\begin{proposition}
If $P$ is a polyhedral cone in general position, then there is $\delta'\equiv\delta_P \mod 2\pi$ such that $(\sigma_P, \delta')\in\Tn$. Conversly, if $(\sigma,\delta)\in \Tn$, then there is a polyhedral cone $P$ in general position, such that $\sigma_P=\sigma$ and $\delta_P\equiv \delta\mod 2\pi$.
\end{proposition}

\begin{proof}
The first implication is true by construction. The converse holds by construction of $\Tn$ and by Proposition~\ref{thm:sphericaltriangle}.
\end{proof}

\begin{proposition}\label{thm:tTn}
	$\tTn$ is a real analytic submanifold of $\Omega_n$ of dimension $2n-3$.
\end{proposition}
 
\begin{proof}
	The function $\gtn$ is smooth. Denote $r_1^i$, $r_2^i$ and $r_3^i$ the rows of $\D \gtn$ corresponding to the $i$-th copy of $\D g_3$ in $D \gtn$ and $s_1,\ldots,s_n$ the last $n$ rows of $\D \gtn$. For $i\in\set{1,\ldots,n-2}$, and $j\in\set{1,\ldots,n}$ let $\lambda_1^i$, $\lambda_2^i$, $\lambda_3^i$ and $\mu_j$ be elements in $\R$ such that
\begin{align}
	\sum_{i=1}^{n-2} \left(\lambda_1^i r_1^i + \lambda_2^i r_2^i + \lambda_3^ir_3^i\right)
	+\sum_{j=1}^n \mu_j s_j &= 0.\label{eq:lambdarimusj}
\end{align}
Note that for every $j\in\set{1,\ldots, n}$ the column $\pby{\gtn}{\delta_j}$ contains only one nonzero entry, which is in the row $s_j$ and is equal to $1$. Therefore, for every $j\in\set{1,\ldots,n}$ Equation~\ref{eq:lambdarimusj} implies
\begin{align}
	\mu_j=0\label{eq:mujzero}
\end{align}
Furthermore, on $\Omega_n$ the nonzero entries of the column $\pby{\gtn}{\alpha_i}$ are in the rows $r_1^i$ and $s_1$ for each $i\in\set{1,\ldots,n-2}$. Therefore, by Equation~\ref{eq:lambdarimusj} and Equation~\ref{eq:mujzero}, we get that for all $i\in\set{1,\ldots,n-2}$ 
\[
	\lambda_1^i = 0.
	\]
By the analogous argument for the columns $\pby{\gtn}{\beta_i}$ and $\pby{\gtn}{\gamma_i}$, we get for every $i\in\set{1,\ldots,n-2}$ 
\[
	\lambda_2^i=\lambda_3^i=0.
	\]
Hence, $\D \gtn$ is of full rank $3(n-2) +n$ on $\Omega_n$ and in particular $0$ is a regular value of $\gtn$. By the preimage theorem it follows that $\tTn$ is a real analytic submanifold of $\Omega_n$ of dimension
\[
	\dim \tTn = \dim\Omega_n - \rk \D\gtn = 2n-3.
\]
\end{proof}

\begin{remark}\label{rmk:eliminateDgn}
\end{remark}

\begin{proposition}\label{thm:tTn=Tn}
$\tTn$ and $\Tn$ are diffeomorphic. Furthermore, $\Tn$ is real analytic.
\end{proposition}

\begin{proof}
Note that the projection
\begin{align*}
	p_n\colon\tTn &\longrightarrow \Tn\\
(\sigma,\delta,\alpha,\beta,\gamma,c)&\longmapsto (\sigma,\delta)
\end{align*}
is smooth and surjective. It remains to construct a smooth inverse function. Let therefore $(\sigma,\delta)\in\Tn$. Then, by definition of $\Tn$ there are 
\[
	(\alpha,\beta,\gamma,c)\in p(\Omega_n),
\]
where $p$ is the projection onto the corresponding factors of $\Omega_n$ and such that
\begin{align}
	\gtn(\sigma,\delta,\alpha,\beta,\gamma,c) = 0\label{eq:defcond}.
\end{align}
We show that $\alpha$, $\beta$, $\gamma$ and $c$ are smooth fuctions of $\sigma$ and $\delta$ by induction on their components $\alpha_i$, $\beta_i$, $\gamma_i$ and $c_i$.

For $i=1$ note that Equation~\ref{eq:defcond} contains the equation
\begin{align}
	\cos c_1 =\cos\sigma_1\cos\sigma_2 + \sin\sigma_1\sin\sigma_2\cos\delta_2.\label{eq:cosc1}
\end{align}
Since the cosine function is invertible on $(0,\pi)$ Equation~\ref{eq:cosc1} implies that $c_1$ is a smooth function $c_1(\sigma,\delta)$ of $\sigma$ and $\delta$. Furthermore, note that by Equation~\ref{eq:defcond} $\gamma_1=\delta_2$ is a smooth function $\gamma_1(\sigma,\delta)$ of $\sigma$ and $\delta$. In addition Equation~\ref{eq:defcond} contains the equation
\begin{align}
	\cos\sigma_2 = \cos\sigma_1\cos c_1 + \sin\sigma_1\sin c_1\cos\alpha_1.\label{eq:alpha1}
\end{align}
Since $\delta_2=\gamma_1$ is given, we know whether $(\alpha_1,\beta_1,\gamma_1)\in (0,\pi)^3$ or whether $(\alpha_1, \beta_1,\gamma_1)\in (\pi,2\pi)^3$. The cosine function is invertible both on $(0,\pi)$ and on $(\pi,2\pi)$. Therefore, since $c_1=c_1(\sigma,\delta)$ is a smooth function of $\sigma$ and $\delta$ Equation~\ref{eq:alpha1} implies that $\alpha_1$ is a smooth function $\alpha_1(\sigma,\delta)$ of $\sigma$ and $\delta$. By an analogous argument $\beta_1$ is a smooth function $\beta_1(\sigma,\delta)$ of $\sigma$ and $\delta$.

Suppose now that $\alpha_i$, $\beta_i$, $\gamma_i$ and $c_i$ are smooth functions of $\sigma$ and $\delta$ for $i\in\set{1,\ldots, m-1}$ and $1\leq m\leq n-2$. Then Equation~\ref{eq:defcond} contains the equation $\gamma_m = \delta_{m+1} -\beta_{m-1}$, which defines by induction hypothesis $\gamma_m$ as a smooth function $\gamma_m(\sigma,\delta)$ of $\sigma$ and $\delta$. Furthermore, by an argument identical to the argument for the case $m=1$ but with shifted indices, we find that $c_m$, $\alpha_m$ and $\beta_m$ are smooth functions $\alpha_m(\sigma,\delta)$, $\alpha_m(\sigma,\delta)$ and $\beta_m(\sigma,\delta)$ of $\sigma$ and $\delta$.

Hence, we have shown, that $\alpha$, $\beta$, $\gamma$ and $c$ are smooth functions $\alpha(\sigma,\delta)$, $\beta(\sigma,\delta)$, $\gamma(\sigma,\delta)$ and $c(\sigma,\delta)$ of $\sigma$ and $\delta$. Define
\begin{align*}
	\psi_n\colon \Tn&\longrightarrow \tTn\\
	(\sigma,\delta)&\longmapsto (\sigma,\delta,\alpha(\sigma,\delta),\beta(\sigma,\delta),\gamma(\sigma,\delta),c(\sigma,\delta)).
\end{align*}
Then $\psi_n$ is a smooth inverse function of $p_n$, so that $\psi_n$ and $p_n$ are diffeomorphisms.
\end{proof}

\subsection{The collection of polyhedral cones of a polyhedron}\label{sec:The_collection_of_polyhedral_cones}

Let $K$ be a combinatoric and $k\in K_0$  be a zero-simplex. We note that a triangulation of a polyhedral cone realizing $C(k)$, can be determined independent of the geometric realization. Let $\tau$ therefore be a choice of a triangulation for every zero-simplex $k\in K_0$. Then we denote $\alpha^k,\beta^k,\delta^k,c^k,\sigma^k$ and $\delta^k$ the variables for angles $\alpha_i^k,\beta_i^k,\delta_i^k$ of the spherical triangles given by $\tau$ at $k$, $c_i^k$ be the variables for lengths of the diagonals of $\tau$ at $k$, and $\sigma_i^k$ and $\delta_i^k$ be the variables for the edgelengths and angles of the spherical $n_k$-gon to be parametrized. Let furthermore
\[
	\Omega_{k,\tau}:= (\Omega_{n_k})^{K_0},
\]and denote $\alpha,\beta,\gamma,c,\sigma$  and $\delta$ the vectors containing $\alpha^k,\beta^k,\gamma^k,c^k\sigma^k$ and $\delta^k$ for every zero-simplex $k\in K_0$. Then, we can define with slight abuse of notation the following.
\begin{definition}
Let $K$ be a combinatoric and $\tau$ be a choice of triangulation for every $k\in K_0$. Let
	\[
	\gtle\colon \Omega_{K,\tau} \longrightarrow \vartimes_{k\in K_0} \R^{3(n_k-2)+ n_k}
	\]
be defined by
\[
	\gtle(\sigma,\delta,\alpha,\beta,\gamma,c)(\set{k})
	:= \widetilde{g}_{n_k}\left(\sigma^k,\delta^k,\alpha^k,\beta^k,\gamma^k,c^k\right),
	\]
and set
\[
	\tTle := \left(\gtle\right)^{-1}(0).
\] 
\end{definition}

\begin{proposition}\label{thm:tTle}
Let $K$ be a combinatoric and $\tau$ be a triangulation of $K_0$. Let $g$ be the genus of the surface $\abs{K}$. Then,
\[
	(\sigma,\delta,\alpha,\beta,\gamma,c)\in \tTle
	\Leftrightarrow \forall k\in K_0\,\,\,
	(\sigma^k,\delta^k,\alpha^k,\beta^k,\gamma^k,c^k) \in \widetilde{\mathcal{T}}_{n_k}.
\] Furthermore, $\tTle$ is a smooth submanifold of $\Omega_{K,\tau}$ of dimension
\[
	 \dim \tTle = 3F + 6g - 6 = 2E + 6g -6.
	\]
\end{proposition}

To prove this proposition we will eliminate the linearized system of equations $\D\gtle =0$. The major challenge in doing so will be to decide on an order in which we eliminate the variables and which equation we use at what moment. In terms of the matrix $\D\gtle$ that means finding a good ordering of its columns and rows. Once a good ordering is found the proof will follow from the structure of the reordered matrix $\D\gtle$. Note that $\D\gtle$ consists of copies of $\D\gtn$, for different $n=n_k$ depending on $k\in K_0$, and recall from Remark~\ref{rmk:eliminateDgn} that we can eliminate $\D\gtn$ in terms of any subset of cardinality $2n-3$ of the partial derivatives by surface and by dihedral angles. Putting it differently, we can determine all auxiliary variables and any set of three of the set of dihedral angles and surface angles in the linear system $\D\gtn =0$ in terms of the remaining surface and dihedral angles. Hence, we can eliminate $\D\gtle$ by eliminating for every $k\in K_0$ three of the variables $\sigma^k$ and $\delta^k$ and all auxiliary variables $\alpha^k$, $\beta^k$, $\gamma^k$ and $c^k$. However, since every dihedral angle occurs as the dihedral angle of two polyhedral cones, and we want to eliminate every copy of $\D\widetilde{g}_{n_k}$ separately from the others, we need to make sure that no $\delta_{ij}$ gets chosen to be eliminated twice.  Equivalently, to eliminate $\D\gtle$ we need to choose for every $k\in K_0$ three elements of the set
\[
	\Set{\text{ One-simplices of $K$ incident in $k$ }}
	\cup
	\Set{\text{ Corners of $K$ at $k$ }},
\]
in such a way that our choices  for each $k\in K_0$ are disjoint. To do so we will dualize $K$ below and we will make our ``choices'' by giving ``colors'' to the dual object $K^\ast$.

\begin{definition}
Let $K$ be a combinatoric. The \mayhem{dual graph} of $K$ is the abstract simplicial one-complex $K^\ast$, such that  
\begin{enumerate}
\item there is a bijection $\ast\colon K_2 \longrightarrow K_0^\ast$, and 
\item two zero-simplices of $K^\ast$ are adjacent in $K^\ast$ if and only if their preimages under $\ast$ are adjacent in $K$ as two-simplices of $K$.
\end{enumerate}
We denote $\set{ijk}^\ast := \ast(\set{ijk})$ and for $k\in K_0$ we call the $n_k$-circuit $C(k)^\ast$, by severe abuse of concepts, a \mayhem{two-simplex of $K^\ast$} and denote it by $k^\ast$. We denote the set of two-simplices of $K^\ast$ by $K^\ast_2$.
\end{definition}

\begin{remark}
The dual graph $K^\ast$ of a combinatoric $K$ is well defined up to adjacency preserving bijections.
\end{remark}

\begin{definition}Let $K$ be a combinatoric. A \mayhem{path of two-simplices} of $K^\ast$ is an ordered set of two-simplices of $K^\ast$ such that any two consecutive elements contain a common one-simplex.
\end{definition}

\begin{remark}
Let $K$ be a combinatoric. Every one-simplex of $K^\ast$ is of the form
\[
	\Set{\set{ijk}^\ast,\set{jkl}^\ast}
\] and is contained in exactly two two-simplices, namely $j^\ast$ and $k^\ast$ of $K^\ast$, which we henceforth call \mayhem{adjacent}. Every two-simplex of $K^\ast$ contains at least three zero-simplices of $K^\ast$. Furthermore, $K^\ast$ is a trivalent connected graph and for any two one-simplices $s_1$ and $s_2$ of $K^\ast$ there is a path $\gamma$ of two-simplices of $K^\ast$ connecting them, that is, such that the first element of $\gamma$ contains $s_1$ and the last element of $\gamma$ contains $s_2$.
\end{remark}

\begin{definition}
Let $K$ be a combinatoric and $k_1^\ast$ and $k_2^\ast$ be two two-simplices of $K^\ast$ that are adjacent in a one-simplex $s$ . Then an \mayhem{arrow} on $s$ is a bijective map $a\in\set{k_1,k_2}^{\set{0,1}}$. We say that $a$ \mayhem{points inward to $a(1)$} and that it \mayhem{points outward of $a(0)$}.
\end{definition}

\begin{definition}
Let $K$ be a combinatoric. Then, the set of \mayhem{corners} of $K^\ast$ is defined by
\[
	C_{K^\ast}:=\left\{(ijk)\middle|	
			\begin{array}{c}
	(ijk)\text{ is an order of }\set{i,j,k}\subset K_0^\ast, \\
	\text{ where } \set{ij},\set{jk}\in K_1^\ast
			\end{array}\right\}\Big/ \sim,
\]
where $\sim$ identifies an order $(ijk)$ with the order $(ijk)$. We say that $(ijk)\in C_{K^\ast}$ is a \mayhem{corner of $K^\ast$ at $j$}. Furthermore, for a two-simplex $k^\ast$ of $K^\ast$, we call the set of corners $(ijl)$ of $K^\ast$ where $\set{ij},\set{jl}\in k^\ast$ the \mayhem{corners of $k^\ast$}.
\end{definition}

\begin{definition} Let $K$ be a combinatoric. A \mayhem{coloring} of $K^\ast$ is a subset $C\subset C_{K^\ast}$ together with a set $A$ of arrows on the one-simplices of $K^\ast$. We call the elements of $C\cup A$ the \mayhem{colors}, call the corners contained in $C$, and the one-simplices contained in $A$ \mayhem{colored}, and say that a one-simplex of $K^\ast$ is \mayhem{colored by} the arrows on it in $A$. 
\end{definition}

\begin{definition}
A coloring $(C,A)$ of $K^\ast$ is called \mayhem{admissible} if for every two-simplex $k^\ast$ of $K^\ast$ it contains exactly three elements of the following
\[
\set{ \text{Arrows pointing outward of $k^\ast$}}\cup \set{\text{Corners of $k^\ast$}}.
\]
\end{definition}

Starting with any specific combinatoric $K$ one can easily find many different admissible colorings of $K^\ast$. It is also easy to see that there always is an admissible coloring. For instance, since every two-simplex of $K^\ast$ contains at least three zero-simplices, one can pick three corners of every two-simplex of $K^\ast$. For the purpose of proving Proposition~\ref{thm:tTle} this kind of coloring is already sufficient. However, to prove Theorem~\ref{thm:tT} we will want an admissible coloring that colors as many one-simplices of $K^\ast$ by exactly one arrow, which is what following conjecture would achieve.

\begin{conjecture}[Elimination Pattern]\label{thm:elimination}
Let $K$ be a combinatoric, and $g$ be the genus of $\abs{K}$. Then there is an admissible coloring of $K^\ast$ such that
\begin{enumerate}
\item If $g=0$, then there is exactly one arrow on every one-simplex of $K^\ast$ and exactly six corners of $K^\ast$ are colored. Furthermore, for any two adjacent zero-simplices $i$ and $j$ of $K^\ast$ the colored corners can be chosen as the union of the corners of $i$ and the corners of $j$.
\item If $g\geq 1$, then there is exactly one arrow on all except $6g-6$ of the one-simplices of $K^\ast$ and no corners are colored.
\end{enumerate}
\end{conjecture}

\begin{proof}[Proof of Proposition~\ref{thm:tTle} modulo the Elimination Pattern Conjecture]
Let 
\[
l\colon\set{1,\ldots,V}\rightarrow K_0
\] be a bijective map and $(C,A)$ be an admissible coloring of $K^\ast$ as obtained by the conclusion of the Elimimation Pattern Conjecture~\ref{thm:elimination}. For $k\in K_0$ denote $C_k$ the set of colored corners of $k^\ast$, and $A_k$ the set of one-simplices of $k^\ast$ that is colored with an outward pointing arrow. Let $\mathcal{C}_n$ be the set of columns of $\D\gtle$ that consists of the partial derivatives
\[
	\pby{\gtle}{\alpha^k}, \pby{\gtle}{\beta^k}, \pby{\gtle}{\gamma^k},\pby{\gtle}{c^k}, \pby{\gtle}{\sigma_{ikj}},\text{ and }\pby{\gtle}{\delta_{ik}}
\]
for $k=l(n)$, all $(ikj)^\ast \in C_k$ and all $\set{jk}^\ast \in A_k$. Denote $\mathcal{R}_n$ the set of rows of $\D\gtle$ of the components $\gtle(\set{l(n)})$ of $\gtle$. Denote
\[
	\mathcal{C}_\text{rest} := \set{\pby{\gtle}{\sigma},\pby{\gtle}{\delta}}\setminus \cup_{\ell=1}^V \mathcal{C}_\ell.
\]
 Choose an ordering of the columns of $\D\gtle$ such that
\[
	\D\gtle =\Parens{\mathcal{C}_1\mid \ldots \mid \mathcal{C}_V\mid \mathcal{C}_\text{rest}}
\]
and such that for $\ell\in\set{1,\ldots,V}$ the columns
\[
	\mathcal{D}_\ell :=\set{\pby{\gtle}{\delta}} \cap \mathcal{C}_\ell
\]
appear after the other columns in $\mathcal{C}_\ell$ in $\D\gtle$. Rearrange the rows of $\D\gtle$ as follows
\[
	\D\gtle = \begin{pmatrix} \mathcal{R}_1\\ \vdots \\ \mathcal{R}_V\end{pmatrix}.
\]
With the chosen ordering of columns and rows of $\D\gtle$ we henceforth consider $\mathcal{C}_\ell$ and $\mathcal{R}_\ell$ as submatrices of $\D\gtle$. 

Then, by the definition of $\gtle$ the submatrix $\mathcal{R}_1$ contains only non-zero entries in the columns that are the partial derivatives by the variables in the vectors $\alpha^k$, $\beta^k$, $\gamma^k$, $c^k$, $\sigma^k$ and $\delta^k$, where $k=l(1)$ and their intersection with $\mathcal{R}_1$ is $\D\gtn$, where $n=n_k$ is the valency of $k$. Hence, by Remark~\ref{rmk:eliminateDgn} we can fully eliminate $\mathcal{R}_1$ in such a way that the block $\mathcal{R}_1 \cap \mathcal{C}_1$ becomes upper triangular. Denote $\mathcal{C}_1'$, respectively
 $\mathcal{D}_1'$, the columns obtained from $\mathcal{C}_1$, respectively $\mathcal{D}_1'$ by this operation and likewise $\mathcal{R}_1'$ the rows obtained from $\mathcal{R}_1$. Then $\mathcal{C}_1'\setminus\mathcal{D}_1'$ is upper triangular. Furthermore, $\mathcal{D}_1'$ contains a quadratic upper triangular submatrix with the same number of columns as $\mathcal{D}_1$ and that is contained in rows of $\D\gtle$ that contain otherwise only zeros. Hence, the other nonzero entries of $\mathcal{D}_1'$ can be eliminated using row operations of $\D\gtle$ without changing $\D\gtle \setminus \mathcal{C}_1'$. Then the submatrix $\mathcal{C}_1''$ obtained from $\mathcal{C}'_1$ by this operation is upper triangular. We successively repeat the argument for $\mathcal{C}_2$, \ldots ,$\mathcal{C}_V$ and thereby fully eliminate $\D\gtle$. 

We compute the rank of $\D\gtle$ by counting its independent columns $\cup_{\ell=1}^{V} \mathcal{C}_\ell$ as follows
\[
\rk\D\gtle = \sum_{k\in K_0} \Parens{3(n_k-2) +(n_k -3)} + 3V,
\]
where the first term counts the number partial derivatives in the variables $\alpha^k$, $\beta^k$, $\gamma^k$ and $c^k$ and the second term counts the cardinality of $\cup_{\ell=1}^V \mathcal{D}_\ell$, which is the same as the total number of colors of $(C,A)$. Similarly, we find
\[
\dim\Omega_{K,\tau}^\text{le} =  E + 3F + \sum_{k\in K_0} \Parens{3(n_k-2) +(n_k -3)},
\] 
and thus by the Preimage Theorem, confer e.g.~\cite{guillemin_pollack:differential_topology:2010}, we compute
\begin{align}
\dim \tTle      &= \dim \Omega_{K,\tau}^\text{le} -\rk \D\gtle\\
		&= 3F + E -3V\\
		&= 3F + 6g -6 	
\end{align}
where the last equation holds by Euler's Polyhedron Formula applied to $V$ and because in any abstract simplicial two-complex $K$ we have that $3F =2E$.
\end{proof}

\section{Wrapping it up}\label{sec:Wrapping_it_up}

	In this section we reconstruct a simply connected polyhedron from its polyhedral cones and its faces. We define the moduli-space $\tT$ of geometric realizations up to similarities of $\R^3$ of a simply connected polyhedron with combinatoric $K$ correspondingly. We recompute its dimension depending on the combinatorc $K$ modulo the Elimination Pattern Conjecture~\ref{thm:elimination}, and reparametrize $\tT$ as $\T$ in terms of dihedral and face angles only.


\subsection{Geometric realizations of simply connected polyhedra}\label{subsec:Geomteric_realizations_of_simply_connected_polyhedra_locally_in_general_position}

To parametrize a simply connected polyhedron we parametrize its collection of polyhedral cones using $\gtle$ and its collection of faces using $\gtin$. Below we show that this is sufficient. 

\begin{definition}
Define
\[
	\gt\colon \Omega_K \times\R_+^{K_1} 
	\longrightarrow 
	 \left(\R^3 \right)^{K_2} \vartimes_{k\in K_0}\R^{3(n_k-2) + n_k}\]
by
\[
	\gt(\sigma,\delta,\alpha,\beta,\gamma,c,x)(\set{k})
	:=\begin{pmatrix}
		\gtin (\sigma,x)\\
		\gtle(\sigma,\delta,\alpha,\beta,\gamma,c)(\set{k})
	\end{pmatrix},
	\]
and set
\[
	\tT := \left( \gt \right)^{-1}(0)
	\quad\text{ and }\quad
	\T := p(\tT),
	\]
where $p$ is the projection onto the $\sigma$- and $\delta$-factors.
\end{definition}

\begin{remark} Note that, by definition $\tT = \tTle \cap \tTin$.
\end{remark}
\begin{vartheorem}[Main Theorem]\label{thm:main}
Let $P\colon\abs{K}\to \R^3$ be a polyhedron that is locally in general position, then there is $\delta'\equiv \delta_P \mod 2\pi$ such that $(\sigma_P,\delta')\in \T$. If $K$ is a combinatoric such that $\abs{K}$ is simply connected then the converse holds, i.e., if $(\sigma,\delta)\in \T$, then there is a polyhedron $P\colon \abs{K}\to\R^3$ such that $\sigma_P = \sigma$ and $\delta_P \equiv \delta \mod 2\pi$.
\end{vartheorem}

To construct a polyhedron $P\colon \abs{K}\rightarrow \R^3$ with given surface angles $\sigma$ and dihedral angles $\delta$ with $(\sigma,\delta)\in \T$ we will find a nice sequence of zero-simplices $k_1,\ldots,k_n\in K_0$ and construct the polyhedral pieces $P'_i$ that will be the faces of $P$ around its vertex $P(\abs{k_i})$. At the same time we will successively fit them together in a sequence of polyhedral pieces
\[
	P_{i+1} := P'_i \cup P_i,
\]
where by abuse of notation we have identified all maps with their image in $\R^3$. We stop the process after finitely many steps at $P= P_N$ for some $N\in \N$. 

To ensure that we can take the union of $P'_i$ and $P_i$ in such a way that they define a well defined map $P_{i+1}$, we need to choose the sequence $k_1,\ldots k_n\in K_0$ in such a way that the union of the cones
\[
	D(k_1,\ldots, k_n) := \bigcup_{i=1}^n \abs{C(k_i)}
\]
has no handles. We achieve this by means of the following two lemmas.

\begin{lemma}\label{thm:adding_the_right_cone_at_the_right_moment}
Let $K$ be a combinatoric such that $\abs{K}$ is simply connected and $k_1,\ldots,k_n\in K_0$ be zero-simplices such that $ D(k_1,\ldots,k_n)$ is a disc and $\abs{K}\setminus D(k_1,\ldots,k_n)\neq\varnothing$. Then there is a zero-simplex $k_{n+1}\in\partial D(k_1,\ldots,k_n)$ such that $D(k_1,\ldots k_{n+1})$ has one connected component and is contractible.
\end{lemma}

\begin{lemma}\label{thm:disc_with_marked_bdry}
Let $D$ be a disc and $S=\partial D$ be marked with points $x_1,\ldots,x_n\in S$. Then, there is no collection $\Gamma$ of paths in $D$ starting and ending in $\set{x_1,\ldots,x_n}$ such that
\begin{enumerate}
\item Every marked point is starting or ending point for some path in $\Gamma$ and for any two distinct points $x_i,x_j\in\set{x_1,\ldots,x_n}$ there is at most one path $\gamma\in\Gamma$ starting in $x_i$ and ending in $x_j$ or vice versa.
\item No two paths of $\Gamma$ cross.
\item Every path $\gamma\in\Gamma$ starts and ends in distinct points.  
\item No path $\gamma\in\Gamma$ starts and ends in $x_i,x_j\in\set{x_1,\ldots,x_n}$ such that $x_i$ and $x_j$ are adjacent in $S$.
\end{enumerate}
\end{lemma}

\begin{proof}
We prove this lemma by induction on $n$. For $n=1$ any collection of paths in $D$ starting and ending in $\set{x_1}$ violates condition $(3)$. Suppose now the statement of the lemma is true for any disc with $m\in\set{1,\ldots,n}$ marked points on its boundary. Let $D$ be a disc with $n+1$ marked points $x_1,\ldots,x_{n+1}\in S=\partial D$. Let $\Gamma$ be a collection of paths in $D$ starting and ending in $\set{x_1,\ldots,x_{n+1}}$ that satisfies conditions (1)-(3). We show that then $\Gamma$ violates condition $(4)$. Let $\gamma\in\Gamma$. Without loss of generality it starts in $x_1$ and ends in $x_i$ for $i\in\set{2,\ldots,n+1}$. The union $S'$ of $\gamma$ with the segments of $S$ from $x_1$ to $x_2$, from $x_2$ to $x_3$, \ldots, and from $x_{i-1}$ to $x_i$ is a simple closed curve and by the Jordan Curve Theorem bound a disc $D'\subset D$. Then $S'=\partial D'$ is marked by $i<n+1$ points $x_1, \ldots, x_i$. Let $\Gamma'\subset \Gamma$ be the subset of paths that start or end in $\set{x_1,\ldots,x_i}$. If $i=2$, then $\gamma$ starts in $x_1$ and ends in $x_2$ that are adjacent so that $\Gamma$ violates condition $(4)$. If $i\neq 2$, then by condition $(1)$ there are paths of  $\Gamma$ starting in $\set{x_2,\ldots,x_i}$ and therefore $\Gamma'\neq \varnothing$. Since no two paths of $\Gamma$ cross by condition $(2)$ every path that starts or ends in $D'$ does not cross $\gamma$ and therefore by the Jordan Curve Theorem stays contained in $D'$. It follows that the paths in $\Gamma'$ all start and end in $\set{x_1,\ldots,x_n}$ and that we can apply the induction hypothesis to $\Gamma'$, $D'$ and $S'$. Since $\Gamma'$ satisfies conditions $(1)$-$(3)$ as a subset of $\Gamma$ it follows that $\Gamma'$ violates condition $(4)$. This means that there is either a path $\gamma'\in \Gamma'$ that starts and ends in $\set{x_1,x_i}$, or that starts and ends in $\set{x_j,x_{j+1}}$ for $j\in\set{2,i-1}$. The first case implies that $\gamma$ and $\gamma'$ are two distinct paths in $\Gamma$ from $x_1$ to $x_i$, which is impossible because $\Gamma$ satisfies condition $(1)$. In the second case $\gamma'$ is a path in $\Gamma$ between two adjacent markings of $S$ so that $\Gamma$ violates condition $(4)$.
\end{proof}

\begin{proof}[Proof of Lemma~\ref{thm:adding_the_right_cone_at_the_right_moment}]
Since $D(k_1,\ldots,k_n)$ is a closed disc in the sphere $\abs{K}$, its complement
\[
	D:=\abs{K}\setminus D(k_1,\ldots,k_n)
\]
is an open disc. Its boundary $S:=\partial D$ is an edgepath of $\abs{K}$ whose vertices mark finitely many points on $S$. If the lemma is not true we will find below in the remainder of the proof for every vertex $\abs{k}$ of $S$ an edge that starts in $\abs{k}$ and ends in a vertex $\abs{i}$ of $S$ that is not adjacent to $k$ in $S$. Since edges are paths that do not cross each other this is impossible by Lemma~\ref{thm:disc_with_marked_bdry}.

First note that, for every vertex $\abs{k}\subset S$,
\[
	\abs{k}\in D(k_1,\ldots,k_n) \cap \abs{C(k)}.
\]
Since $D(k_1,\ldots,k_n)$ and $\abs{C(k)}$ are connected, it follows that their union $D(k_1,\ldots,k_n,k)$ is connected. Therefore, if the statement of the lemma is not true, then for every vertex $\abs{k}\subset S$ there is a closed path $\gamma_k$ that is non-contractible in $D(k_1,\ldots,k_n,k)$. Furthermore, $\gamma_k$ is homotopic in $D(k_1,\ldots,k_n,k)$ to a closed simple non-contractible edgepath, which we also denote by $\gamma_k$. The intersection $D\cap \gamma_k$ is nonempty, since otherwise $\gamma_k\subset D(k_1,\ldots,k_n)$, and since $D(k_1,\ldots,k_n)$ is simply connected by hypothesis, $\gamma_k$ would be contractible in $D(k_1,\ldots,k_n,k)$. Furthermore, since $\gamma_k$ is a simple edgepath its connected components in $\overline{D}$ are edgepaths starting and ending at pairwise distinct vertices of $S$.

Denote by $\gamma^{ij}_k$ the connected component of $\gamma_k$ in $\overline{D}$ starting in the vertex $\abs{i}$ of $S$ and ending in the vertex $\abs{j}$ of $S$. There is a connected component $\gamma_k^{ij}$ of $\gamma_k$ in $\overline{D}$ such that either $\abs{i}$ or $\abs{j}$ are not adjacent to $\abs{k}$ in $S$; for otherwise, since $\gamma_k^{ij}\subset \abs{C(k)}$ and $\abs{C(k)}$ is simply connected, $\gamma_k^{ij}$ would be homotopic in $\abs{C(k)}$ to $\abs{ik}\cup\abs{kj}$ and $\gamma_k$ contractible in $D(k_1,\ldots,k_n,k)$.

Since every connected component $\gamma_k^{ij}$ of $\gamma_k$ in $\overline{D}$ is a subset of $\abs{C(k)}$ and $\abs{C(k)}$ is simply connected there is furthermore a path $\widetilde{\gamma}_k^{ij}$ that is homotopic to $\gamma_k^{ij}$ in $\abs{C(k)}$ and that contains $\abs{k}$. Let $\widetilde{\gamma}_k$ be the path homotopic to $\gamma_k$ in $D(k_1,\ldots,k_n,k)$ obtained by replacing $\gamma_k^{ij}$ by $\widetilde{\gamma}_k^{ij}$. 

Consider the connected components $\gamma_k^{ij}$ of $\gamma_k$ in $\overline{D}$ such that $\abs{i}$ or $\abs{j}$ is not adjacent to $\abs{k}$ in $S$. Without loss of generality it is $\abs{i}$ that is not adjacent to $\abs{k}$ in $S$. Consider the segments $\overline{\gamma}_k^{ij}$ of $\widetilde{\gamma}_k^{ij}$ starting in $\abs{k}$ and ending in $\abs{i}$. Then, since $i\in\St k$ there is an edge $\abs{ik}\subset \abs{C(k)}$ and since both $\abs{ik}$ and $\overline{\gamma}_k^{ij}$ are paths in $\abs{C(k)}$ that start in $\abs{k}$ and end in $\abs{i}$ they are homotopic in $\abs{C(k)}$. 

Let $\overline{\gamma_k}$ be the path obtained from $\widetilde{\gamma_k}$ by replacing the segments $\overline{\gamma}_k^{ij}$ of $\widetilde{\gamma}_k^{ij}$ by the edges $\abs{ik}$ whenever $\abs{i}$ is not adjacent to $\abs{k}$ in $S$. If for all connected components $\gamma_k^{ij}$ of $\gamma_k$ for which $\abs{i}$ is not adjacent to $\abs{k}$ in $S$, the edge $\abs{ik}$ was contained in $D(k_1,\ldots,k_n)$ already, then $\overline{\gamma}_k$ would be contractible in $D(k_1,\ldots,k_n,k)$. But $\gamma_k$ is homotopic to $\overline{\gamma}_k$ in $D(k_1,\ldots,k_n)\cup \abs{C(k)}$ and therefore $\gamma_k$ would be contractible in $D(k_1,\ldots,k_n)\cup \abs{C(k)}$. It follows that, for every vertex $\abs{k}$ of $S$ there is an edge $\abs{ik}$ in $\overline{D}$ and such that $\abs{i}$ is not adjacent to $\abs{k}$ in $S$.
\end{proof}

\begin{lemma}\label{thm:intersection_of_new_cone_with_Dn}
Let $K$ be a combinatoric such that $\abs{K}$ is simply connected and let $k_1,\ldots, k_N\in K_0$ be zero-simplices of $K$ such that for every $n\in\set{1,\ldots,N-1}$ the set $D(k_1,\ldots,k_n)$ satisfies the hypothesis of the previous Lemma~\ref{thm:adding_the_right_cone_at_the_right_moment}, and $k_{n+1}\in K_0$ satisfies the conclusions of the previous Lemma~\ref{thm:adding_the_right_cone_at_the_right_moment}. Then,
\[
	\abs{C(k_{n+1})} \cap D(k_1,\ldots,k_n) = \cup_{l=1}^m f_l \cup \gamma_1\cup \gamma_2,
\]
where $f_1,\ldots,f_m$ are consecutively adjacent faces of $\abs{K}$ with the common vertex $\abs{k_{n+1}}$ and $\gamma_1$ and $\gamma_2$ are two edge-paths of $\abs{K}$ that start at a vertex of $f_1$, respectively $f_m$, not contained in $f_2,\ldots,f_{m-1}$. That is, the above intersection looks e.g. as follows,
\[
	\op{conelooks}{\put(40,15){$f_1$}
			\put(57,15){$f_2$}
			\put(70,21){$f_3$}
			\put(2,2.5){$D(k_1,\ldots,k_n)$}
			\put(20,47){$\color{cyan}\gamma_1\color{black}$}
			\put(78,36){$\color{cyan}\gamma_2\color{black}$}
			\put(50,72){$\abs{K}\setminus D(k_1,\ldots,k_n)$}
			}
\]
where the depicted cone is $\abs{C(k_{n+1})}$, the shaded region is $D(k_1,\ldots, k_n)$ and the unshaded region is its complement in $\abs{K}$.
\end{lemma}

\begin{proof}
We first note that $D(k_1) = \abs{C(k_1)}$ is such that for every zero-simplex $\abs{s}\in\partial D(k_1)$ in its boundary, the faces it shares with $\abs{C(s)}$ are two adjacent faces, and that $D(k_1)\setminus \abs{C(s)}$ is edge-path-connected and $\abs{C(s)}\cap D(k_1)$ satisfies the conclusion of the lemma. Now suppose that for some $n\in\set{1,\ldots,N-1}$ the set $D(k_1,\ldots, k_n)$ is such that for every zero-simplex $s\in \partial D(k_1,\ldots,k_n)$ the faces it shares with $\abs{C(s)}$ are consecutively adjacent, and that the closure of $D(k_1,\ldots, k_n)\setminus \abs{C(s)}$ is edge-path-connected.

Then, if $x_1, x_2\in\partial\abs{C(k_{n+1})}\cap D(k_1,\ldots,k_n)$ are two zero-simplices, there is a simple edge-path $\gamma'$ in the closure of $D(k_1,\ldots, k_n)\setminus \abs{C(k_{n+1})}$  from $x_1$ to $x_2$. Let $\beta_1$ and $\beta_2$ be the two distinct edgepaths from $x_1$ to $x_2$ on $\partial \abs{C(k_{n+1})}$. Then $\gamma:=\gamma' \cup \beta_1$ is a simple closed path and by the Jordan Curve Theorem bounds two disjoint open discs $D_1$ and $D_2$ on $\abs{K}$ such that $D_1\cup\gamma\cup D_2= \abs{K}$. Since $\beta_2$ intersects $\gamma$ only in its starting and ending points, its interior is contained in either $D_1$ or in $D_2$.  Let without loss of generality $D_2$ be the disc containing the interior of $\beta_2$. Then $D_2$ also contains the interior of $\abs{C(k_{n+1})}$. It follows, that $D_1\cap \abs{C(k_{n+1})} =\varnothing$.

 Since by the previous Lemma~\ref{thm:adding_the_right_cone_at_the_right_moment} the set $D(k_1,\ldots,k_{n+1})$ is contractible, it follows that $D_1$ or $D_2$ is contained in $D(k_1,\ldots, k_{n+1})$. If $D_1\subset D(k_1, \ldots, k_{n+1})$, then, since $D_1\cap\abs{C(k_{n+1})}=\varnothing$ it follows that $D_1\subset D(k_1,\ldots,k_n)$ and therefore the closure $D_1\cup\gamma$ of $D_1$ is also contained in $D(k_1,\ldots,k_n)$. Since $\beta_1\subset \gamma$ we have that $\beta_1 \subset D(k_1,\ldots,k_n)$. If $D_1$ is not contained in $D(k_1,\ldots,k_{n+1})$, then $D_2$ is contained in $D(k_1,\ldots, k_{n+1})$. However, then $\abs{C(k_{n+1})}$ is a homotopy in $D_2$ from $\beta_1$ to $\beta_2$ and therefore $D_2$ can be retracted to the complement of the interior of $\abs{C(k_{n+1})}$. Hence , $\gamma'\cup \beta_2$ bounds a disc in $D(k_1,\ldots,k_n)$ and in particular we have that $\beta_2\subset D(k_1,\ldots,k_n)$.

It follows, that the intersection $D(k_1,\ldots,k_{n+1})\cap \abs{C(k_{n+1})}$ looks as claimed. Furthermore, for any zero-simplex $\abs{s}$ of $\partial D(k_1,\ldots, k_{n+1})$, the faces in the intersection of $\abs{C(s)}$ with $D(k_1,\ldots, k_{n+1})$ are again consecutively adjacent, and $D(k_1,\ldots, k_{n+1})\setminus \abs{C(s)}$ is edge-path connected. This completes the inductive argument and proves the lemma.

\end{proof}

\begin{proof}[Proof of Main Theorem~\ref{thm:main}]
The first part of the statement follows by construction. 

For the second part, let $K$ be a combinatoric such that $\abs{K}$ is simply connected and let $(\sigma,\delta)\in \T$. By the definition of $\T$ the vector $\sigma$ has as entries the angles of a flat-cone surface $\ell\colon K_1\to \R_+$  triangulated by $K$. By Lemma~\ref{thm:intrinsicdiffeo}, $\ell$ is determined uniquely from $\sigma$ and the choice of one of its edgelength. Choose a one-simplex $\set{ij}\in K_1$ and set the length of $\set{ij}$ equal to $\ell_{ij}=1$. Then, Lemma~\ref{thm:intrinsicdiffeo} determines $\ell$ from $\sigma$ and $\ell_{ij}=1$ by $\varphi_{ij}(\sigma,1) =(\sigma,\ell)$. In the following we construct a polyhedron $P$ with surface angles $\sigma$, dihedral angles $\delta_P \equiv \delta \mod 2\pi$ and edgelengths $\ell$.

Let $k_1=i$. We construct a simplexwise affine and simplexwise injective map $P_1\colon \abs{C(k_1)}\rightarrow \R^3$ as follows. By the definition of $\T$ and by Proposition~\ref{thm:tTn} there is a polyhedral cone $C_1$ with underlying abstract simplicial cone $C(k_1)$ and surface angles $\sigma^{k_1}$ and dihedral angles $\delta_{C_1}\equiv\delta^k \mod 2\pi$. Send $\abs{k_1}$ under $P_1$ to the tip of $C_1$. Furthermore, send $\abs{k}$ for every $k\in \St{k_1}$ under $P_1$ to the unique point on the edge $C_1(\abs{k_1 k})$ of $C_1$ of distance $\ell_{k_1 k}$ from $P_1(\abs{k_1})$. We have constructed the image of $P_1$ of all the zero-simplices of $\abs{C(k_1)}$. Define $P_1$ on any two-simplex $\abs{s k_1 t}$ of $\abs{C(k_1)}$ as the affine map defined by the images under $P_1$ of its zero-simplices $\abs{s}$, $\abs{k_1}$ and $\abs{t}$. This also defines the image of $P_1$ of the one-simplices of $\abs{C(k_1)}$. By the Side-Angle-Side Congruence Theorem for euclidean triangles the surface angles of $P_1$ are indeed the components of $\sigma^{k_1}$ and its edgelengths are given by the restriction of $\ell$ to the one-simplices of $C(k_1)$. By constuction the dihedral angles of $P_1$ are $\delta_{P_1}\equiv\delta^k \mod 2\pi$. By construction $P_1$ is simpexwise affine and since the entries of $\ell$ are positive it is simplexwise injective.

By Lemma~\ref{thm:adding_the_right_cone_at_the_right_moment} there is a subset $\set{k_1,\ldots,k_N}\subset K_0$ such that for each $n\in\set{1,\ldots,N}$ the subsets $\abs{D_n}:=D(k_1,\ldots,k_n)$ of $\abs{K}$ are simply connected, such that $\abs{k_{n+1}}\in \partial \abs{D_n}$ and such that $\abs{D_N}=\abs{K}$. Let now $n\in\set{1,\dots, N-1}$ and suppose there is a simplexwise affine and simplexwise injective map $P_n\colon \abs{D_n}\rightarrow \R^3$ whose surface angles are given by the restriction of $\sigma$ to the corners of $D_n:=\cup_{i=1}^n C(k_i)$, whose dihedral angles coincide up to $2\pi$ with $\delta$ on the one-simplices of $D_n$ and whose edgelengths are given by $\ell$ restricted to the one-simplices of $D_n$. 

In the same way as we contructed $P_1$ we construct a simplexwise affine and simplexwise injective map $P_n'\colon \abs{C(k_{n+1})}\rightarrow \R^3$ whose surface angles are given by $\sigma$ restricted to the corners of $\abs{C(k_{n+1})}$, whose dihedral angles are $\delta$ modulo $2\pi$ restricted to the one-simplices of $\abs{C(k_{n+1})}$ and whose edgelengths are given by $\ell$ restricted to the one-simplices of $\abs{C(k_{n+1})}$. By Lemma~\ref{thm:intersection_of_new_cone_with_Dn} the intersection of $\abs{C(k_{n+1})}$ with $D(k_1,\ldots,k_n)$ is given by
\[
	\abs{C(k_{n+1})} \cap D(k_1,\ldots,k_n) = \cup_{l=1}^m f_l \cup \gamma_1\cup \gamma_2,
\]
where $f_1,\ldots, f_m$ are consecutively adjacent faces and $\gamma_1$ and $\gamma_2$ are egdepaths starting at a vertex of $f_1$, and of $f_m$, respectively.
The image of $f_1\cup\ldots \cup f_m$ under $P_n$ is a set of consecutively adjacent triangles in $\R^3$ with the same edgelengths and dihedral angles beween them, as the image of $f_1\cup\ldots\cup f_m$ under $P_n'$. It follows, that $P_n(f_1\cup\ldots\cup f_m)$ and $P_n'(f_1\cup\ldots\cup f_m)$ are congruent in  $\R^3$. Let $T$ be the euclidean isometry that takes the latter to the former and from now on also denote by $P_n'$ the map $T\circ P_n'$.

Then, by construction $P_n$ and $P_n'$ coincide on $f_1\cup\ldots\cup f_m$. We now argue, that $P_n$ and $P_n'$ also coincide on $\gamma_1$ and $\gamma_2$. If $\gamma_1$ and $\gamma_2$ are the empty set, or if they are contained in $f_1\cup \ldots \cup f_m$, then the images of $P_n$ and $P_n'$ coincide on $\gamma_1$ and $\gamma_2$. Otherwise, consider the vertex $\abs{k}$ of $f_1$ in which $\gamma_1$ starts suppose that $\abs{l_1}$ is the adjacent vertex to $\abs{k}$ on $\gamma_1$ that is not contained in $f_1\cup\ldots \cup f_m$. Then the union of the image $P_n(\abs{C(k)}\cap D(k_1,\ldots,k_n))$ and the image $P_n'(\abs{C(k)}\cap\abs{C(k_{n+1})})$ is a set of consecutively adjacent triangles in $\R^3$ starting in the line segment $P_n(\abs{kl_1})$ and ending in the line segment $P_n'(\abs{kl_1})$.  It is uniquely determined by the edge $P_n(\abs{k\ell_1})$, the angles of the triangles and the dihedral angles between any two consecutive triangles. However, the faces of
\[
 C:=P_n(\abs{C(k)}\cap D(k_1,\ldots,k_n)) \cup P_n'(\abs{C(k)}\cap\abs{C(k_{n+1})})
\]
 have by construction the surface angles and dihedral angles beween adjacent faces of a polyhedral cone. Therefore, $C$ defines a polyhedral cone and its first and last edges $P_n(\abs{kl_1})$ and $P_n'(\abs{kl_1})$ coincide. We repeat the argument for the next edge of $\gamma_1$ and sucessively for the remaining edges of $\gamma_1$ and obtain that $P_n$ and $P_n'$ coincide on $\gamma_1$. The same argument works for $\gamma_2$.

Define $P_{n+1}\colon D(k_1,\ldots, k_{n+1})\rightarrow \R^3$ as $P_n$ on $D(k_1,\ldots, k_n)$ and as $P_n'$ on $\abs{C(k_{n+1})}$. By the previous arguments the map $P_{n+1}$ is well defined. By construction it is simplexwise affine, simplexwise injective, and has surface and dihedral angles given by the maps $\sigma$ restricted to the corners of $D(k_1,\ldots,k_{n+1})$, and $\delta$ restricted to $D(k_1,\ldots, k_{n+1})$ and componentwise taken modulo $2\pi$.
\end{proof}

\begin{remark}\label{rmk:redundancy}
Note that in the above construction the last two faces are implied as soon as the surrounding ones have been constructed. In particular, the equations $\gtin(\sigma,x)(\set{ijk}) = 0$ for the two corresponding two-simplices $\set{ijk}$ are redundant, and may be dropped from the definition of $\tT$.
\end{remark}

\subsection{Dimension of the manifold of geometric realizations}\label{subsec:Dimension_of_the_manifold_of_geometric_realizations}

In the following we indicate how to recompute the dimension of $\T$ modulo the Elimination Pattern Conjecture~\ref{thm:elimination}.

\begin{theorem}[modulo EPC~\ref{thm:elimination}]\label{thm:tT}
Let $K$ be a combinatoric such that $\abs{K}$ is simply connected. Then $\tT$ is an analytic submanifold of $\Omega_{K,\tau}\times\R_+^{K_1}$ of dimension
\[
	\dim \tT = E.	
\]
\end{theorem}
\begin{proof}[Proof modulo EPC~\ref{thm:elimination}]
By definition $\tT$ is the locus of $\Omega_{K,\tau}\times\R_+^{K_1}$ where
\begin{align}
\gtin(\sigma,x) = 0 \quad \text{ and }\quad \gtle(\sigma,\delta,\alpha,\beta,\gamma,c) = 0\label{eq:locus}.
\end{align}
By  Remark~\ref{rmk:redundancy}, there are two adjacent two-simplices $\set{ijk}$ and $\set{jkl}\in K_2$, such that if we remove the equations
\[
	\gtin(\sigma,x)(\set{ijk})=0 \quad \text{ and }\quad
	\gtin(\sigma,x)(\set{jkl})=0
\]
from the system of equations~\ref{eq:locus}, then the new system of equations still defines the same set of solutions $\tT$. Define a new function
\[
	f_K^\text{in}\colon (0,\pi)^{C_K}\times\R_+^{K_1} \longrightarrow (\R^3)^{K_2\setminus\set{\set{ijk},\set{jkl}}}
\]
by
\[
	f_K^\text{in}(\sigma,x)(\set{k}) :=\gtin(\sigma,x)(\set{k}),
\]
and a function
\[
	f_K\colon \Omega_{K,\tau}\times\R_+^{K_1} \longrightarrow (\R^3)^{K_2}\vartimes_{k\in K_0} \R^{3(n_k -2) +n_k}
\]
by
\[
	f_K(\sigma,\delta,\alpha,\beta,\gamma,c,x)(\set{k}) 
	:= \begin{pmatrix} f_K^\text{in}(\sigma,x)\\\gtle(\sigma,\delta,\alpha.\beta,\gamma,c)(\set{k})\end{pmatrix}.
\]
Then, by the above $\tT= f_K^{-1}(0)$. By the Elimination Pattern Conjecture~\ref{thm:elimination} there is an admissible coloring $(C,A)$ of $K^\ast$ that satisfies the conclusions of the Lemma~\ref{thm:elimination} and such that its colored corners $C$ are the corners of $\set{ijk}^\ast$, together with the corners of $\set{jkl}^\ast$. By Proposition~\ref{thm:tTle} and by its proof $\D\gtle$ is of full rank and furthermore it has a maximal submatrix $M$ of full rank so that the columns of $M$ contained in $\pby{f_K}{\sigma}$ are the partial derivatives $\pby{f_K^\text{in}}{\sigma_{rst}}$ where $(rst)$ are the corners of either $\set{ijk}$ or of $\set{jkl}$. It folows that $M$ and $\D f_K^\text{in}$ as submatrices of $\D f_K$ share no columns of $\D f_K$. Since by Proposition~\ref{thm:intrinsic} we have that $\D f_K^\text{in}$ is of full rank on $(f_K^\text{in})^{-1}(0)$, we conclude that $\D f_K$ is of full rank on $f_K^{-1}(0)$.

By the Preimage Theorem~\cite{guillemin_pollack:differential_topology:2010} it follows that $\tT$ is a smooth submanifold of $\Omega_{K,\tau}\times \R_+^{K_1}$of dimension
\[
	\dim\tT = \dim \left(\Omega_{K,\tau}\times\R_+^{K_1}\right) - \rk \D f_K.
\]
We compute
\[
	\dim\left(\Omega_K\times\R_+^{K_1}\right) = \sum_{k \in K_0} \left(3(n_k -2) +(n_k-3)\right) + 3F + 2E,
\]
and
\begin{align*}
	\rk \D f_K 	&=\rk \D f_K^\text{in} + \rk \D\gtle \\ 
		   	&=(\rk \D\gtin -6) + \rk\D\gtle \\ 
			&= 3F -6 + \left(\sum_{k \in K_0} \left(3(n_k -2) +(n_k-3)\right) + 3V\right)
\end{align*}
and obtain
\begin{align}
	\dim\tT &= 3F + 2E - (3F +3V -6)\\
		&= (E -3V) + E + 6\\
		&= - 6 + E +6\label{eq:sixg}\\
		&= E,
\end{align}
where Equation~\ref{eq:sixg} holds by Euler's Polyhedron Formula applied to $V$, and because in every simplicial two-complex $K$ we have that $2E=3F$.
\end{proof}

\begin{theorem}\label{thm:diffeoTK} [modulo EPC~\ref{thm:elimination}]
Let $K$ be a combinatoric. For every $\set{ij}\in K_1$ there is a diffeomorphism
\[
	\psi_{ij} \colon \T \times \R_+ \longrightarrow \tT. 
	\]
In particular, if $\abs{K}$ is simply connected, then $\T$ is a real analytic submanifold of $(0,\pi)^{C_K}\times \R_+^{K_1}$ of dimension
\[
	\dim \mathcal{T}_{K,\tau} = E - 1.
	\]
\end{theorem}

\begin{proof}[Proof modulo EPC~\ref{thm:elimination}]
The projection map
\begin{align*}
	\pi_{ij}\colon \tT &\longrightarrow \T\times\R_+\\
	(\sigma,\delta,\alpha,\beta,\gamma,c,x)&\longmapsto (\sigma,\delta,x_{ij})
\end{align*}
is smooth. Let $(\sigma,\delta,\ell)\in\T\times\R_+$. Then, for every $k\in K_0$ there are $(\alpha_k,\beta_k,\gamma_k,c_k)$ in an appropriately projected domain, such that 
\[
(\sigma^k,\delta^k,\alpha_k,\beta_k,\gamma_k,c_k)\in\widetilde{\mathcal{T}}_{n_k}.
\] Recall from the proof of Proposition~\ref{thm:tTn=Tn} that $\alpha_k$, $\beta_k$, $\gamma_k$ and $c_k$ can be expressed as smooth functions of $\sigma^k$ and $\delta^k$. Hence $\alpha$, $\beta$, $\gamma$ and $c$ can be expressed as smooth functions $\alpha(\sigma,\delta)$, $\beta(\sigma,\delta)$, $\gamma(\sigma,\delta)$ and $c(\sigma,\delta)$ of $\sigma$ and $\delta$. Then the function
\begin{align*}
	\psi_{ij}\colon \T\times\R_+ &\longrightarrow \tT\\
	(\sigma,\delta,\ell) &\longmapsto (\sigma,\delta,\alpha(\sigma,\delta),\beta(\sigma,\delta),\gamma(\sigma,\delta),c(\sigma,\delta),d_{ij}(\sigma,\ell)),
\end{align*}
where $d_{ij}$ is the function defined in the proof of Proposition~\ref{thm:intrinsicdiffeo}, defines a smooth inverse for $\pi_{ij}$.
\end{proof}

\bibliography{flexie.bib}{}
\bibliographystyle{plain}
\end{document}